\documentclass[11pt]{amsart}
\hoffset= - 0.75 in

\usepackage{amsfonts, amsmath, amssymb, amsthm}
\usepackage{latexsym}
\usepackage{graphics}
\newtheorem*{defn}{Definition}
\newtheorem{thm}{Theorem}
\newtheorem{qn}{Question}

\newtheorem{obs}{Observation}
\newtheorem{lem}{Lemma}
\newtheorem{cor}{Corollary}

\newcommand{\ZZ}{\mathbb{Z}}
\newcommand{\RR}{\mathbb{R}}

\newcommand{\NN}{\mathbb{N}}

\begin{document}

\title{Cayley Compactifications of Abelian Groups}
\author{Mike Develin}
\address{Department of Mathematics, UC-Berkeley, Berkeley, CA 94720-3840}
\curraddr{2706B Martin Luther King Jr. Way, Berkeley, CA 94703}
\date{\today}
\email{develin@math.berkeley.edu}

\begin{abstract}

Following work of Rieffel~\cite{rieffel}, in this document we define the Cayley compactification of a
discrete group $G$ together with a set of generators $S$. We use algebraic methods in the general case to
construct an explicit presentation of Cayley compactifications. In the particular case of $\ZZ^m$, we use
geometric methods to demonstrate a strong connection between the Cayley compactification and the polytope
$P^\Delta$ polar to $P=\text{conv}(S)$. With their explicit description, Cayley compactifications thus 
form a nice class of compactifications of these objects.

\end{abstract}

\maketitle

\section{Introduction}

In~\cite{rieffel}, Rieffel discusses the structure of group
$C^\star$-algebras as compact topological spaces. Among other things, he
constructs an embedding of the discrete group $\ZZ^n$, given a set of
generators, into a compact space by means of adding points at infinity.
This construction generalizes to arbitrary groups and generators with
weights, and has interesting algebraic and geometric structure. In this
paper, we take a closer look at this structure.

First, using geometric methods, we consider the case of $\ZZ^n$. We 
establish a topological correspondence between the Cayley compactification 
and the polytope polar to $P=\text{conv}(S)$, where $S$ is the generating 
set in question; a simple consequence of this is that the Cayley 
compactification is in fact compact. We give an essentially complete 
description of the behavior of the points that are added at infinity to 
compactify $\ZZ^n$.

Moving to the case of an arbitrary group, we utilize algebraic methods, 
specifically that of standard pair decomposition, to give a general 
description of the Cayley compactification of an arbitrary group and 
generating set, provided that the cost vector giving the weights of the 
generators is generic. 

We begin by giving the definitions of the objects in question, along with 
a few relatively straightforward lemmas. In Section~\ref{examples} we 
compute a number of examples which illustrate the behavior of Cayley 
compactifications. Section~\ref{boundary} and Section~\ref{grobner} deal 
with the geometric and algebraic approaches respectively.

\section{Definitions and useful lemmas}

Given a group $G$ (all groups in this paper will be abelian) and a finite set of
generators $S=S^{-1}$, we can define the \textbf{Cayley graph} $(G,S)$ to have vertex set equal to the
elements of $G$, with two vertices $u$ and $v$ being connected by an edge if $uv^{-1}\in S$. The 
generators may be given with weights, in which case the length of this edge is defined to be the weight 
of the generator $uv^{-1}$. 

This construction induces a metric space structure on $G$, using the obvious distance function $d(u,v)$ of
distance between the corresponding vertices in the Cayley graph. We can then define a 
class
of functions $\phi_{y,z}$, where $y$ and $z$ are vertices of the graph, via $\phi_{y,z}(x)=d(x,y)-d(x,z)$. We
will be dealing only with Cayley graphs, but this definition makes sense for any graph $G$. For a graph $G$,
define $\mathbf{B_G}$ to be the $\ZZ$-subalgebra of all functions $f:V(G)\rightarrow \ZZ$ generated by
constant functions and the functions $\phi_{y,z}$.  $B_G$ then encodes the metric space structure of the 
graph.

We state a pair of trivial lemmas concerning the $\phi_{y,z}$.

\begin{lem}\label{trivial}
$\phi_{w,y}+\phi_{y,z}=\phi_{w,z}$; $\phi_{z,y}=-\phi_{y,z}$.
\end{lem}

\begin{lem}\label{genonly}

If we have a Cayley graph $(G,S)$, then the $\phi_{y,z}$ are generated by elements of the form
$\phi_{x,sx}$ for $s\in S$.  Indeed, if $T$ is any set of generators for $G$, the $\phi_{y,z}$ are
contained in the $\ZZ$-module generated by elements of the form $\phi_{x,tx}$ for $t\in T$.

\end{lem}

We now define the central concept which allows us to construct the Cayley compactification, that of a geodesic.

\begin{defn}
Suppose $G$ is a graph, and we have a map $\phi:D\rightarrow G$, $D\subseteq \NN$, such that for all $e<f$, $e,f\in D$, 
$d_G(\phi(e),\phi(f)) = f-e$. Then we say that $\gamma = \text{im}(\phi)\subset G$ is a \textbf{geodesic} in $G$. If $D$ is 
an unbounded subset of $\NN$, we call $\gamma$ an \textbf{infinite geodesic}; if $D=\NN$, we call $\gamma$ a \textbf{full 
geodesic}. $\gamma$ is a full geodesic if and only if $\gamma$ consists of the vertices of an infinite path with the 
following condition: for each $u,v$ in the path, the subsegment of the path from $u$ to $v$ provides a shortest path in $G$ 
between the two vertices. 
\end{defn}

We then have the following definition and lemma, which explains why geodesics are useful at 
all.

\begin{defn}
For any infinite geodesic $\gamma = \{\gamma_0,\gamma_1,\ldots\}$, we define the operator $v_\gamma$ via 
$v_\gamma(f) = \lim_{n\rightarrow\infty}f(\gamma_n)$.
\end{defn}

It is not a priori clear why this limit exists, but it does so in the case we are concerned with, as we show in 
the following 
lemma.

\begin{lem}
For $f\in B_G$, and for $\gamma = \{\gamma_0,\gamma_1,\ldots\}$ any infinite geodesic, $v_\gamma(f)$ exists.
\end{lem}
\begin{proof}
It suffices to show it for the generators of $B_G$, since $v_\gamma$ is obviously an ring homomorphism; it's trivial 
for the constant functions. By Lemma~\ref{genonly}, it suffices to show it for functions of the form 
$\phi_{y,\gamma_0}$.

Now, because of the triangle inequality, we know that $-d(y,\gamma_0)\le \phi_{y,\gamma_0}(\gamma_i)\le d(y,\gamma_0)$. We 
claim that the function $\phi_{y,\gamma_0}(\gamma_i)$ is nonincreasing as $i$ increases, which will finish the proof since 
its value is an integer in a bounded range, and hence must converge if it is nonincreasing. For any $i<j$, 
we have $\phi_{y,\gamma_0}(\gamma_j) = d(\gamma_j,y) - d(\gamma_j,\gamma_0)$. But by the geodesic property, we know that 
$d(\gamma_j,\gamma_0) = d(\gamma_j,\gamma_i)+d(\gamma_i,\gamma_0)$; by the triangle inequality, we also have 
$d(\gamma_j,y)\le d(\gamma_j,\gamma_i) + d(\gamma_i,y)$. Plugging this back in, we get $\phi_{y,\gamma_0}(\gamma_j) \le 
d(\gamma_i,y)-d(\gamma_i,\gamma_0) = \phi_{y,\gamma_0}(\gamma_i)$ as desired, which completes the proof.
\end{proof}

So, for each $\gamma$, we have a valuation (which is an algebra homomorphism) $v_\gamma:B_G\rightarrow \ZZ$. In addition, for 
each $g\in G$, we have the natural valuation $g^\star:B_G\rightarrow \ZZ$ given by $g^\star(f)=f(g)$. This brings us to our 
main definition.

\begin{defn}
Let $(G,S)$ be a Cayley graph. Then the \textbf{Cayley compactification} $C(G,S)$ consists of the set of valuations 
$\{v_\gamma\}\cup \{g^\star\}$, where $\gamma$ ranges over all infinite geodesics and $g$ ranges over all elements of $G$. 
The topology is the $\text{weak}^\star$ topology, where a sequence of valuations $v_1,v_2,\ldots$ converges to $v$ if for any 
function $f\in B_G$, $v_1(f),v_2(f),\ldots$ converges to $v(f)$. Note that if $\gamma=\{\gamma_1,\ldots\}$ is a geodesic, 
then by definition $\gamma_1^\star,\gamma_2^\star,\ldots$ converge to $v_\gamma$.
\end{defn}

Indeed, this definition makes sense for any graph, not just Cayley graphs. However, Cayley graphs, especially those of 
integer lattices $\ZZ^n$, provide an extremely tractable subclass of infinite graphs with which to work. 

Justifying the name, Cayley compactifications are indeed actually compact.

\begin{thm}
For any group $G$ and generating set $S$, the Cayley compactification $C(G,S)$ is compact. 
\end{thm}

We will give two proofs of this theorem via the two main techniques presented in this paper. In
Section~\ref{boundary} we wil prove it for the case of $\ZZ^n$ with unweighted generators using geometric
methods; in Section~\ref{grobner} we will give an algebraic proof which works in the case of a general
group $G$ with arbitrarily weighted generators.

Finding geodesics of Cayley graphs is usually quite easy. The tricky part is determining what the
equivalences are; this amounts to determining when two geodesics induce the same valuation $v_\gamma$,
as the following two lemmas dispatch the other cases.

\begin{lem}
If $g$ and $h$ are two distinct elements of $G$, $g^\star$ and $h^\star$ are distinct valuations, and thus $g$ and $h$ 
correspond to different points in the Cayley compactification.
\end{lem}
\begin{proof}
Plugging into the definition, it is easy to see that $g^\star(\phi_{g,h})=-h^\star(\phi_{g,h})\neq 0$.
\end{proof}

\begin{lem}
Suppose $g$ is an element of $G$, and $\gamma$ is a geodesic. Then $g^\star$ and $v_\gamma$ are distinct valuations.
\end{lem}

\begin{proof}

From the construction before, where we showed that $v_\gamma$ converges on any $\phi$, we can find
sufficiently large $i$ so that for all $j>i$, there exists a shortest path from $\gamma_j$ to $g$ going
through $\gamma_i$. (We can assume $g\neq \gamma_i$, since if so we can just pick a larger $i$.) But now we
have $g^\star(\phi_{g,\gamma_i}) = -\gamma_j^\star(\phi_{g,\gamma_i})$ for all $j>i$, and so
$g^\star(\phi_{g,\gamma_i})=-v_\gamma(\phi_{g,\gamma_i})\neq 0$.

\end{proof}

So each element of $G$ corresponds to a distinct point in the Cayley compactification $C(G,S)$, which
justifies the nomenclature of this object as a compactification of the original graph (or group.) The
interesting part comes in determining which geodesics induce distinct valuations, and what the space of
distinct valuations looks like; this varies depending on the set of generators $S$. We define the boundary of
the Cayley compactification, $\partial C(G,S)$, to consist of the valuations $v_\gamma$ inside of the Cayley
compactification. The main question can then be posed as follows:

\begin{qn}
For a given $(G,S)$, what does the boundary $\partial C(G,S)$ look like?
\end{qn}

In general, this is a very difficult question, which we will investigate in Section~\ref{boundary}. We can 
narrow the question down a bit by making the following definition and corresponding observations.

\begin{defn}
If we have two geodesics $\gamma$ and $\eta$ such that $v_\gamma = v_\eta$ as valuations on the function 
algebra $B_G$, then we call $\gamma$ and $\eta$ \textbf{equivalent}. The boundary $\partial C(G,S)$ then consists 
of the equivalence classes of geodesics under this notion of equivalence.
\end{defn}

\begin{obs}
If $\gamma$ and $\eta$ are geodesics that share an infinite set of vertices, then they are equivalent.
\end{obs}

\begin{obs}
When determining the set of $v_\gamma$, it suffices to consider full geodesics, since any infinite geodesic
is contained in a full geodesic simply by taking the union of paths between $\gamma_i$ and $\gamma_j$, where
$i$ and $j$ are consecutive in the indexing set $D$, arranged in increasing order.
\end{obs}

Since the Cayley graphs are translation-invariant, it follows that the natural group action of $G$ on itself 
extends to the Cayley compactification $C(G,S)$. When determining whether or not a sequence of elements 
comprises a full geodesic, the only thing of importance is therefore the differences between successive 
elements; we make the following natural definitions.

\begin{defn}
A geodesic $\gamma$ \textbf{contains} a multiset $T$ of generators if the elements 
$\{\gamma_i-\gamma_{i-1}\}$ 
contain that multiset. The \textbf{direction} of a geodesic $\gamma$ is defined to be the subset of generators 
which appear infinitely often among the elements $\{\gamma_i-\gamma_{i-1}\}$. 
\end{defn}

Still, the space of full geodesics is large, and in general is difficult to parametrize. We can cut down the
space significantly by using lemmas such as the following, which has an easy direct proof. We state this fact
now so as to allow us to compute examples before developing the machinery.

\begin{lem}\label{minimal}
Let $G$ be an abelian group, $S = \{e_1,\ldots,e_n\}$. Then if $e_{i_1}+\cdots+e_{i_k}$ can be written as the
sum of fewer than $k$ generators, no full geodesic $\gamma$ in the Cayley graph $(G,S)$ can contain
$\{e_{i_1},\ldots,e_{i_k}\}$.
\end{lem}

\begin{proof}

Suppose $\gamma$ contains $e_{i_1},\ldots,e_{i_k}$. Pick $i$ sufficiently large that the multiset
$\{\gamma_j-\gamma_{j-1} \mid j\le i\}$ contains all of these elements with multiplicities. By definition of
a geodesic, we know that the shortest path between $\gamma_i$ and $\gamma_0$ goes through all intermediate
$\gamma_j$, having length $i$; in particular, $\gamma_i-\gamma_0$ cannot be expressed as the sum of fewer
than $i$ generators. However, we can replace the $k$ elements $e_{i_1},\ldots,e_{i_k}$ with fewer than $k$
generators in the sum for $\gamma_i-\gamma_0$, a contradiction.

\end{proof}

Having proven some preliminary results which allow us to compute Cayley compactifications in specific cases, we present some examples of Cayley 
compactifications, which illustrate their general structure as well as a variety of behaviors that can occur.

\section{Examples}\label{examples}

\textit{Example 1.} Let $G=\ZZ^2$, $S=\{\pm e_1,\pm e_2\}$ where $e_1 = (1,0)$, $e_2 = (0,1)$ is the standard basis.

If $\gamma$ is a full geodesic, by Lemma~\ref{minimal}, $\gamma$ cannot contain both $e_1$ and $-e_1$ or both $e_2$ and 
$-e_2$. By symmetry, it suffices to consider the case where $\gamma$ contains only $e_1$ and $e_2$ (by reflection we will 
obtain all the other full geodesics.) We break this down into three cases.

1. $\gamma$ contains infinitely many $e_1$ and infinitely many $e_2$. We claim that all such $\gamma$ induce the same 
valuation $v_\gamma$. Because of Lemma~\ref{genonly}, it suffices to check this on functions $\phi_{(x,y),(x+1,y)}$ and 
$\phi_{(x,y),(x,y+1)}$. But for any such $\gamma$ and a given $(x,y)$, eventually all points in $\gamma$ are of the form 
$(a,b)$ where $a>x+1$ and $b>y+1$, and so we will have $d((a,b),(x,y))=(a-x)+(b-y)=\delta$, and similarly $d((a,b),(x+1,y)) = 
(a-(x+1)) + (b-y) = \delta-1$, $d((a,b),(x,y+1) = \delta-1$. So we have $v_\gamma(\phi) = 1$ for both of these cases, 
regardless of what $\gamma$ is in this class. Consequently, all such $v_\gamma$ are the same, and evaluate to $1$ on all such 
functions.

2. $\gamma$ contains infinitely many $e_2$ and finitely many $e_1$, say $k$ of them. Because they necessarily share an 
infinite set of vertices, $\gamma$ is equivalent to the full geodesic $\{(k,z)\mid z\in \NN\}$. Again, we evaluate 
$v_\gamma$ on the functions $\phi_{(x,y),(x+1,y)}$ and $\phi_{(x,y),(x,y+1)}$. The calculation is easy in both cases; we end 
up with $v_\gamma(\phi_{(x,y),(x+1,y)}) = 1$ if $x < k$ and $-1$ otherwise, and $v_\gamma(\phi_{(x,y),(x,y+1)}) = 1$ 
regardless.

3. $\gamma$ contains infinitely many $e_1$ and $k$ copies of $e_2$. This is isomorphic to the previous case, and yields 
$v_\gamma(\phi_{(x,y),(x+1,y)}) = 1$ always, and $v_\gamma(\phi_{(x,y),(x,y+1)}) = 1$ if $y < k$ and $-1$ otherwise.

As you can see, all geodesics obtained in these three cases are distinct; it is 
easy to check that their symmetric images (using $-e_1$ and/or $-e_2$ in place of $e_1$ and $e_2$ respectively) are also distinct from 
each other and from these $v_\gamma$. Figure 1 shows the resulting Cayley compactification.

\begin{figure}\label{figone}
\input{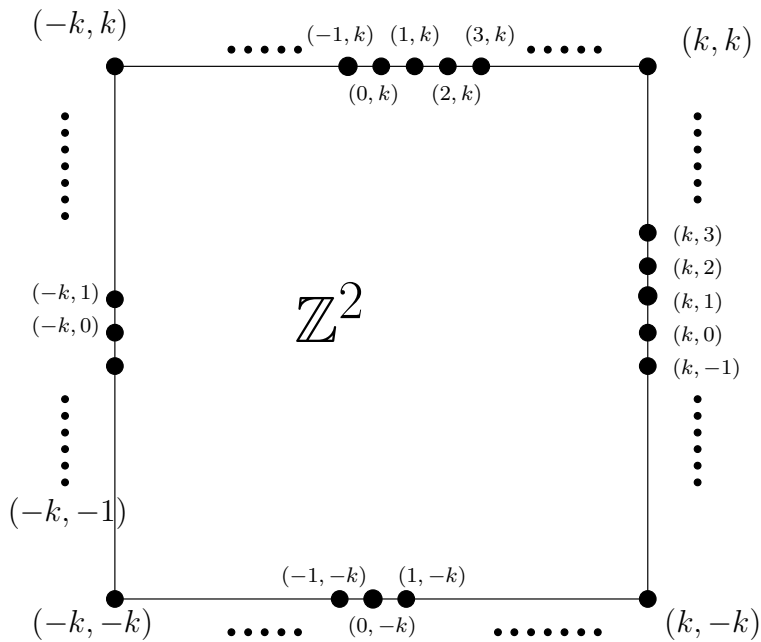}
\caption{The Cayley compactification $C(\mathbb{Z}^2, \{\pm (1,0), \pm
(0,1)\})$. The interior of the square depicted is just an infinite grid,
$\mathbb{Z}^2$ with the discrete topology. The geodesics form a square
compactifying this grid as shown; $(1,k)$ represents the geodesic
consisting of those points, and so on.}
\end{figure}

The drawing of the figure reflects the topology of the situation; the geodesics $\{(k,z)\mid z\in \NN\}$ converge (as $k$ 
gets large) to the geodesic with infinitely many of both $e_1$ and $e_2$.

This example reflects one phenomenon of the Cayley compactification. The boundary $\partial C(G,S)$ is 
the boundary 
of the polytope polar to the one defined by the generators $\pm e_1$ and $\pm e_2$. This is not a coincidence, and it's easy 
to see how it arises; anything "between" $e_1$ and $e_2$ resulted in the same point, while infinitely taking $e_1$ resulted 
in a facet of geodesics, one for each translation of $e_1$ by a finite combination of the other basis vectors.

As the next two examples show, it is not generally the case that the boundary of the Cayley compactification is simply equal to the boundary of 
the polar polytope in some geometric sense. Rather, in general, the boundary of the Cayley compactification has components that look like the 
boundary of the polar polytope, spliced together; we will prove a more rigorous statement to this effect in Section~\ref{boundary}.

\textit{Example 2.} Let $G=\ZZ$, $S = \{\pm 1, \pm 8\}$.

As before, thanks to Lemma~\ref{minimal}, no full geodesic can contain both $\pm 1$, both $\pm 8$, or more than four copies of either $+1$ or
$-1$. Consequently, the only full geodesics contain infinitely many copies of either $+8$ or $-8$. It suffices to discuss the $+8$ case.

If a full geodesic contains finitely many items that aren't $+8$, then as before it is equivalent to some geodesic 
$\gamma_k = \{k+8n\mid n\in\NN\}$, with $k\in \{0,1,\ldots,7\}$. It is straightforward to show that $v_{\gamma_k}\neq 
v_{\gamma_l}$ for $k\neq l$, simply by evaluating them on $\phi_{k,l}$. The geodesics with infinitely many $-8$'s are also 
easily seen to be distinct from these. So, we have eight distinct geodesics in either direction; the 
Cayley compactification is depicted in Figure 2.

\begin{figure}\label{figtwo}
\input{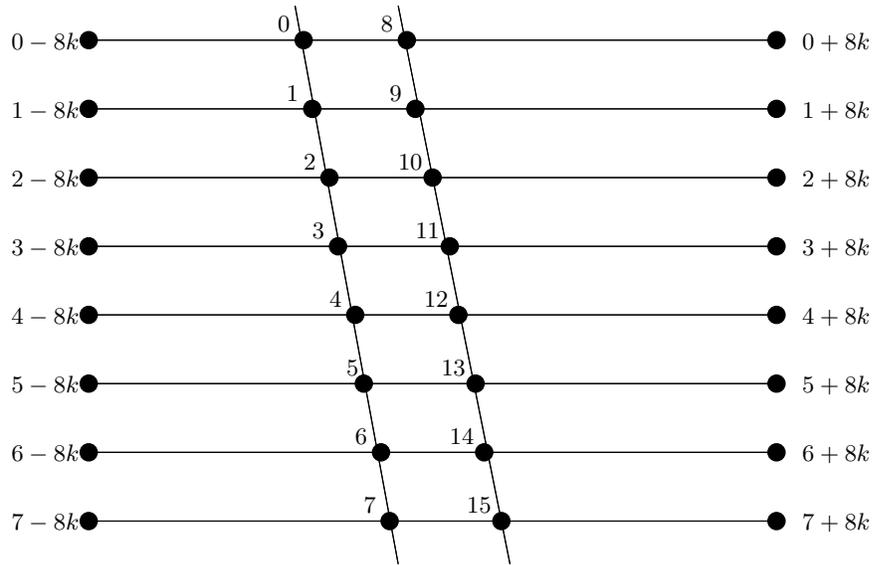}
\caption{The Cayley compactification $C(\mathbb{Z}, \{\pm 1, \pm 8\})$. Again,
the discrete group $\mathbb{Z}$ comprises the interior of the figure, although we
have shown the metric structure by giving a closeup of the Cayley graph.
There
are eight strands, corresponding to each of the eight cosets of
$\mathbb{Z}/8\mathbb{Z}$, each compactified by one geodesic in the $+8$ direction and
one in the $-8$ direction.}
\end{figure}

In this case, we obtain eight copies of the boundary of the polar polytope, which can be represented as $\{k+8t, k-8t\}$ for $0\le k\le 7$. This 
comes from the fact that the lattice generated by $+8$ is not saturated; each of the eight copies corresponds to a coset of this lattice group, 
an element of $\ZZ/8\ZZ$.

\textit{Example 3.} Let $G=\ZZ^2$, $S = \{\pm e_1, \pm e_2, \pm e_3\}$, where we have $e_1 = (1,0)$, $e_2 = (0,1)$, $e_3 = 
(2,2)$.

This example combines both of the previous ones. On the one hand, we have the polar polytope corresponding to the convex hull 
of the basis elements, which is a hexagon. This shows up as follows: the geodesics fall into two categories. One is families; 
you have families in the $(1,0)$ direction, in the $(0,1)$ direction, and in the $(2,2)$ direction, and all of their 
negatives; these are the six 1-faces of the hexagon, corresponding to the four sides of the square from Example 1. The other 
corresponds to the four corners of the square in Example 1; we have six corners now, one corresponding to geodesics with 
infinitely many of $(1,0)$ and $(2,2)$, the adjacent ones corresponding to geodesics with infinitely many of $(2,2)$ and 
$(0,1)$ or $(1,0)$ and $(0,-1)$, and so on. All of these geodesics are equivalent; this can be shown 
directly, and it also falls out as a consequence of results to follow.

On the other hand, as in Example 2, there is also a so-called torsion part in $\partial C(G,S)$, and it shows up in an odd way. Note that we
have two "overlapping" (though of course they don't overlap in the Cayley graph, but in the picture they do) geodesics in the $(2,2)$ direction:
for instance, the geodesics $\{(2x,2x)\}$ and $\{2x+1,2x+1\}$. Indeed, all geodesics in this direction (or its negative) come in such pairs.

As in Example 2, you might think that the boundary would comprise two copies of the polytope (as we had
eight before). However, it is more complicated than that; we do in fact have two copies of each face,
but they fit together in interesting ways. For instance, the subsequences of geodesics with direction
$(2,2)$ which converge to a geodesic containing two generators infinitely often are interlaced, as
depicted in Figure 3.

\begin{figure}\label{figthree}
\input{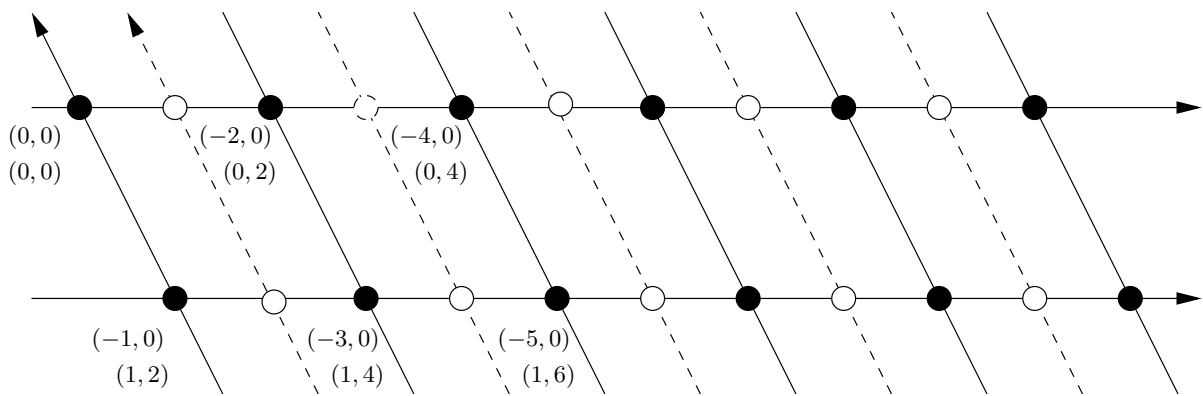}

\caption{Geodesics in the $(2,2)$-direction in the Cayley compactification $C(\mathbb{Z}^2, \{\pm 
(1,0), \pm (0,1), \pm (2,2)\})$. As you can see, there are two limit points in each of two
directions, corresponding to the two geodesics with direction $\{(2,2),(0,1)\}$ and the two with
direction $\{(2,2),(1,0)\}$. The diagonal solid and dashed lines correspond to points converging
to the latter geodesics, while the horizontal lines correspond to points converging to the former
ones. We have given two labels to a number of points; the labels demonstrate the interlacing,
with the top labels showing how these six points are part of a geodesic with direction
$\{(2,2),(1,0)\}$, and the bottom labels showing how they naturally fall into geodesics with
direction $\{(2,2),(0,1)\}$.}

\end{figure}

Having presented these examples which hopefully illustrate the intuitive idea of what Cayley
compactifications look like, we move on to approaches for investigating them. We will present two
approaches: a geometric approach using polytopes, which is simpler and works better in the case of
$\ZZ^m$ (which is embedded in $\RR^m$), and an algebraic approach using standard pairs, which works in a
more general setting.  We first present the geometric approach, which gives all of the results that the
algebraic approach gives, and more, in the case of $\ZZ^m$.

\section{Geometric methods}\label{boundary}
In this section we investigate the shape of the boundary $\partial C(\ZZ^n,S)$ using geometric methods; 
throughout this section $P$ will denote the polytope which is the convex hull of the vectors in the 
generating set $S$. Note that since $S$ generates $\ZZ^n$, 0 is in the interior of $P$.

By truncating a full geodesic $\gamma$, we can assume that the only generators that are contained in
$\gamma$ are contained in its direction, as truncation obviously preserves equivalence class. Therefore,
we can restrict our attention to looking at full geodesics which only contain generators in their
direction; unless otherwise stated, we will assume that this is the case for all geodesics mentioned in
the remainder of this section. Furthermore, we have the following very powerful result.

\begin{thm}\label{basedir}
Suppose two full geodesics $\gamma_1$ and $\gamma_2$ with the same direction $\{t_1,\ldots,t_n\}\subset S$ contain a common point $x$. Then 
$\gamma_1$ and $\gamma_2$ are 
equivalent.
\end{thm}

\begin{proof}

Our method will be to construct a new geodesic $\eta$ which has infinitely many points in common with both $\gamma_1$ and
$\gamma_2$; this geodesic will then be equivalent to both of our original geodesics, and so they must be equivalent to each
other. First of all, we can truncate both $\gamma_1$ and $\gamma_2$ at $x$, and thus assume that $x$ is the first point of both.

To each point $w$ in $\gamma_1$ or $\gamma_2$, we associate an $n$-tuple of nonnegative integers
$(a_1(w),\ldots,a_n(w))$, where $a_i(w)$ is the number of times that the generator $t_i$ is contained in
the subgeodesic stretching from $x$ to $w$ (which is merely a path.) By construction, these $n$-tuples
are nondecreasing in each component as we move along the geodesics in question. Furthermore, as we move
off towards infinity along either $\gamma_j$, for all $i$, $a_i$ will eventually be arbitrarily large,
since each $t_i$ appears infinitely often in $\gamma_j$.

Now, let $\eta = (x_0,y_0,x_1,y_1,\ldots)$, where the $x_i$ are from $\gamma_1$ and the $y_i$ from $\gamma_2$. We select these 
points as follows: $x_0$ is any point from $\gamma_1$. Each $y_i$ is closen so that $a_j(y_i)>a_j(x_i)$ for all $j$, and 
similarly each $x_i$ except the first is chosen so that $a_j(x_i) > a_j(y_{i-1})$ for all $j$. It is then clear that these points 
comprise a geodesic, since applying Lemma~\ref{minimal} to $\gamma_1$ yields that no sum of the $t_i$'s with any multiplicities 
can be expressed as the sum of fewer generators, which is the only obstruction to $\eta$ being a geodesic.
\end{proof}

Theorem~\ref{basedir} justifies the nomenclature in question; a geodesic is completely categorized by
where it starts and which direction it goes in, as one would expect from the geometry of the situation.
Thus, all points which contain a given point $x$ and have the same direction $T$ are equivalent; we will
write such a geodesic as $\gamma = (x,T)$, and refer to $x$ as a \textbf{base} of $\gamma$. The next step 
is to classify possible directions of geodesics;
we start with an easy lemma.

\begin{lem}\label{every}
Let $x$ be a point with rational coordinates in the interior of a rational polytope $P$ (i.e. a polytope 
all of whose vertices have rational coordinates.) Then there exists an affine combination of all the 
vertices of $P$ with positive rational coefficients which is equal to $x$.
\end{lem}
\begin{proof}
First, we note that it suffices to show that for every vertex $v$ of $P$, there exists an affine 
combination of the vertices of $P$ with rational coefficients such that the coefficient of $v$ is 
positive and all other coefficients are nonnegative; we can then just average all of these combinations 
to obtain one where all coefficients are positive.

The rest of the proof is by induction on the dimension of $P$. Take any vertex $v$, and consider the ray
from $v$ to $x$. This ray intersects $P$ beyond $x$, since $x$ is in the interior of $P$; it does so at
some point $y$ in the interior of some proper face $F$, which is of course a polytope itself. By induction, $y$ 
can be expressed as an affine combination of the vertices of $F$ with rational coefficients. But now $x$ 
can be expressed as a rational affine combination of $v$ and $y$, and since $x$ is in the interior of $P$ 
and $y$ is not, the coefficient of $v$ must be positive. We have therefore produced an expression for $x$ 
of the desired form.
\end{proof}

Using this lemma, we can prove a nice theorem about geodesics that classifies them by their directions.

\begin{thm}\label{face}
Let $\gamma$ be a full geodesic. Then $\gamma$ is equivalent to a geodesic whose direction consists of 
the vertices of some proper face $F$ of $P$.
\end{thm}

\begin{proof}
Suppose there are $n$ vertices $v_i$ in the direction of $\gamma$; consider the rational point $x =
(v_1+v_2+\cdots+v_n)/n$. This is in $P$, because $P$ is convex, and we have $nx= v_1+v_2+\cdots+v_n$.  
Because of Theorem~\ref{basedir}, $\gamma$ is equivalent to a geodesic $\eta$ which starts at the 
same point as $\gamma$ and then repeats this block of generators adding up to $nx$ ad infinitum.

Now, $x$ is contained in some face $F$ of $P$, so by Lemma~\ref{every} we can find an affine combination of the vertices of $F$ 
with rational coefficients which is equal to $x$. Multiplying through by the common denominator yields that, for some positive 
integer $l$, we can write $lx$ as a positive integer linear combination of the vertices of $F$, where the coefficients add up to 
$l$.

Consider blocks of $nl$ generators of $\eta$; by construction, each block adds up to $nlx$. Replace each of these blocks with $n$ 
copies of the expression for $lx$ in terms of the vertices of $F$; these copies together will also have $nl$ generators, and 
thus the resulting path $\zeta$ must also be a geodesic, since $\eta$ is. Furthermore, $\zeta$ and $\eta$ share an infinite set 
of vertices, namely those at intervals of $nl$ generators. Consequently, $\zeta$ is equivalent to $\eta$, and since the direction 
of $\zeta$ is clearly the vertex set of $F$, the theorem is nearly proven.

All that remains is to show that no geodesic can have a direction equal to all of $P$. However, since 0 can be written as a 
positive linear combination of all of the vertices of $P$, being in its interior, this is patently impossible.
\end{proof}

Theorem~\ref{basedir} tells us that two geodesics which share a point and which have the same associated face are equivalent; 
however, we have as of yet proved very little about when two geodesics cannot be equivalent. The following theorem is a major 
step in that direction.

\begin{thm}\label{noneq}
Suppose $\gamma$ and $\eta$ are two geodesics with directions equal to the vertex sets of 
distinct faces $F$ and $G$ respectively. Then $\gamma$ and $\eta$ are not equivalent.
\end{thm}

\begin{proof}
Without loss of generality, let $v$ be a vertex in $F$ but not in $G$. Consider the function 
$f = \phi_{\eta_0, \eta_0+v}$; we claim that this function has different limits along $\gamma$ and 
$\eta$.

First, we compute the limit along $\gamma$. For sufficiently large $i$, we have $f(\gamma_i) =
f(\gamma_{i+1})$. Also, as demonstrated in the proof of that fact, for sufficiently large $i$ we have
$d(\gamma_{i+1}, \eta_0) - d(\gamma_{i}, \eta_0) = 1$. Take a pair with $i$ exceeding this bound such
that $\gamma_{i+1}-\gamma_i = v$; we can do this, since $v$ occurs infinitely many times in $\gamma$. We
then have $d(\gamma_{i+1}, \eta_0) = 1 + d(\gamma_{i}, \eta_0) = 1 + d(\gamma_{i}+v, \eta_0+v) = 1+
d(\gamma_{i+1}, \eta_0+v)$, which means that $v_\gamma(f) = 1$.

On the other hand, consider the limit of $f$ along $\eta$. We have $d(\eta_i, \eta_0) = i$ by definition;  
we claim that we always have $d(\eta_i, \eta_0 + v) > i-1$, which completes the proof since $v_\eta(f)$
must then be strictly less than 1. To do this, we consider the face-defining linear functional $\omega$
of $G$, which has $\omega(x) = 1$ for $x\in G$ and $\omega(x)<1$ otherwise.

On the one hand, since every generator contained in $\eta$ has $\omega=1$, we know that $\omega(\eta_i) = 
i + \omega(\eta_0)$. On the other hand, since $v\notin G$, $\omega(\eta_0 + v) < 1 + \omega(\eta_0)$. If 
we had a path of length $i-1$ or less connecting $\eta_0+v$ and $\eta_i$, we would obtain a set of $i-1$ 
generators whose $\omega$-values summed to $\omega(\eta_i)-\omega(\eta_0+v) > i-1$. This is impossible by 
definition of $\omega$, and so the distance between them must be strictly greater than $i-1$, which 
completes the proof.
\end{proof}

We have now shown that for investigating equivalence classes of full geodesics, it suffices to consider directions which are the 
vertex sets of faces of the polytope $P=\text{conv}(S)$. We can reformulate this as a statement about the natural group action of 
$\ZZ^n$ on $C(\ZZ^n, S)$.

\begin{cor}\label{corresp}
The orbits of the natural group action of $\ZZ^n$ on $C(\ZZ^n,S)$ are in one-to-one correspondence with 
the proper faces of the polytope $P = \text{conv}(S)$. 
\end{cor}

\begin{proof}
The orbit of a geodesic is simply the set of geodesics of the same direction. Directions correspond to 
proper faces; it is easy to see that every proper face has a geodesic with that direction, simply by 
considering the face-defining linear functional. The empty set corresponds naturally to the orbit 
consisting of $\ZZ^n$ itself; these "one-point geodesics" have the empty set as their direction.
\end{proof}

Corollary~\ref{corresp} establishes that the components of the boundary of $C(\ZZ^n,S)$ correspond to the components of the
polytope $P=\text{conv}(S)$. However, we can do better: in fact, these components fit together in the same way that the faces of
the polar polytope $P^\Delta$ do. In particular, $P^\Delta$ has the property that the closure (in Euclidean space) of the
interior of the face corresponding to $F$ is precisely the union of the faces $G^\Delta$, where $G$ ranges over all faces
containing $F$ (and thus over all faces for which $G^\Delta\subset F^\Delta$.) $\partial C(\ZZ^n,S)$ has 
the same property.

\begin{thm}\label{topol}
Let $O_F$ be the orbit corresponding to the face $F$ as described in Corollary~\ref{corresp}. Then its closure $\overline{O_F}$
in $C(\ZZ^n,S)$ is equal to $\bigcup O_G$, where the union is taken over all proper faces $G$ of $P$ 
containing $F$.
\end{thm}

\begin{proof}
First, we need to show that for all such $G$, a point from $O_G$ is in $\overline{O_F}$; since the $\ZZ^n$-action is a homeomorphism on 
$C(\ZZ^n,S)$, it then follows that all of $O_G$ is. By abuse of notation, we write $F$ for the direction consisting of the vertices of $F$, and 
similarly for $G$. Let $\{x_i \mid i\in T\}$ be the set of generators in $G$ but not in $F$.

Consider the set of geodesics $(\gamma_0,\gamma_1,\gamma_2,\ldots) = ((0,F), (x^T,F), (2x^T,F),
\ldots)$, where $x^T = \sum_{i\in T} x_i$. We claim that this set of geodesics converges (in the 
topology of 
$C(\ZZ^n,S)$, the weak
topology related to the algebra $B_{G,S}$) to a geodesic with direction $G$. Take any function $f =
\phi_{y,z}$; we need to show that $v_{\gamma_i}(f)$ converges to $v_\eta(f)$.

Let $p_0 = 0$. Then, for each $i$ in increasing order, pick $p_i$ sufficiently large such that
$v_{\gamma_i}(f) = f(p_i)$ and such that for each $x_j\in F$, the multiplicity of $x_j$ in $p_i$ (i.e.  
the number of times $x_j$ occurs in the subgeodesic path from $ix^{T}$ to $p_i$) is greater than the
multplicity of $x_j$ in $p_{i-1}$. Now, it is easy to see by previous results that the $p_i$ form a
geodesic $\eta=(0,G)$, since for each $x_j\in G$, the multiplicity of $x_j$ in $p_i$ is greater than the
multiplicity of $x_j$ in $p_{i-1}$. Since $v_{\gamma_i}(f) = f(p_i)$, the $v_{\gamma_i}(f)$ converge to
$v_\eta(f)$. This is true for every $f$, and so the $v_{\gamma_i}$, which are points in $O_F$, converge
to $\eta\in O_G$ as desired.

To complete the proof, we need to show that if a face $G$ does not contain $F$, then no point in $O_G$ is the limit of a sequence of points in 
$O_F$. As before, if one point is a limit of a sequence of points, then every point must be, so we need only show that there exists a point in $O_G$ 
which is not the limit of a sequence of points in $O_F$. Let $x_i$ be a point in $F$ but not in $G$; let $\omega$ be a linear functional defining 
$G$, so that $\omega(x_j)\le 1$ with equality holding if and only if $x_j$ is in $G$. Note that we have $\omega(x_i) < 1$.

Now, consider the distance difference function $f=\phi_{0,x_i}$, and consider the geodesic $\eta = 
(0,G)$. By the argument we presented in the 
proof of Theorem~\ref{noneq}, we must have $v_\eta(f)\le 0$. On the other hand, for any geodesic $\gamma$ with direction $F$, again as in 
Theorem~\ref{noneq} we have $v_\gamma(f) = 1$, since $x_i$ is contained infinitely many times in $\gamma$. Since convergence is in the weak topology 
corresponding to an algebra containing $f$, $\eta$ cannot be a limit point of elements of $O_F$.

\end{proof}

Theorem~\ref{topol} yields the following corollary, which gives further credence to the appellation of the 
space as a compactification of the original group.

\begin{cor} 
$\ZZ^n = O_\varnothing$ is open and dense in $C(\ZZ^n, S)$.
\end{cor}

\begin{proof}
Density is immediately obvious from Theorem~\ref{topol}. Meanwhile, the complement is simply $\bigcup O_F$, 
where $F$ ranges over all nonempty proper faces of $P$. This is a finite union, so its closure is $\bigcup 
\overline{O_F}$, which looking at each $O_F$ individually and applying Theorem~\ref{topol} is clearly equal to 
$\bigcup O_F$ itself.
\end{proof}

Theorem~\ref{topol} demonstrates that the orbits correspond in this nice topological way to the faces of the 
polar polytope $P^\Delta$. One question remains:  what do the orbits look like? The answer is relatively 
simple.

\begin{thm}
For a face $F$, let $H_F$ be the subgroup of $\ZZ^n$ generated by $T$, the set of points of $S$ lying on $F$. Then the points of the orbit $O_F$ 
are in 
natural bijection with the elements of $\ZZ^n/H_F$.
\end{thm}

\begin{proof}

First of all, note that if a geodesic has direction $F$ and base $B$, then it is equivalent to a geodesic
with direction $T$ and base $B$. This is because for any element $t\in T$, we can write $kt$ as a sum of
vertices of $F$, using $k$ in all. Consequently, we can rewrite that geodesic using infinitely many
copies of $t$. Proceeding for every element of $T$ which is not a vertex of $F$, we obtain the desired
result.

We need to show two things: that two geodesics with direction $T$ are equivalent if their bases differ by an element of $H_F$, and that two 
geodesics with direction $T$ are not equivalent if their bases differ by an element not in $H_F$.

The first claim is easy. If two bases $x_1$ and $y_1$ differ by an element of $H_F$, then we can find $x_2$ and $y_2$, positive
linear combinations of the elements of $T$, such that $x_1+x_2 = y_1+y_2$. We can then shift generators corresponding to $x_2$ to
the front of the geodesic based at $x_1$, and similarly for $y_2$ and $y_1$, to show that both geodesics are equivalent to a 
geodesic with 
base $x_1+x_2$
and direction $T$, and hence equivalent to each other.

Next, suppose that we have two equivalent geodesics $\gamma_1$ and $\gamma_2$ with direction $F$ and bases $x_1$ and $y_1$ respectively. As usual, let
$\omega$ be the face-defining linear functional of $F$, so that we have $\omega(x)\le 1$ for $x\in S$ with equality holding if and only if $x\in T$.
Without loss of generality, by truncating the geodesics we can assume that $x_1$ and $y_1$ satisfy $0\le \omega(y_1) -\omega(x_1) < 1$. First, suppose
$\omega(y_1) > \omega(x_1)$. Then consider the distance difference function $f = \phi_{y_1,x_1}$. By the same argument as in the proof of
Theorem~\ref{noneq}, we must have $v_{\gamma_1}(f) \ge 0$, and $v_{\gamma_2}(f) \le -1$, and so in particular these two cannot be equivalent.

On the other hand, suppose $\omega(y_1) = \omega(x_1)$. Then as before we have $v_{\gamma_1}(f) \ge 0$ and $v_{\gamma_2}(f) \le 0$, and so if the two 
geodesics are equivalent these must in fact both be equal to 0. This means that for sufficiently large $n$, the distance from $(\gamma_1)_n$ to $y_1$ 
is the same as the distance from $(\gamma_1)_n$ to $x_1$. The latter path consists only of elements of $T$ since the direction of $\gamma_1$ is $F$; 
the former path must consist of the same number of elements, and since we have $\omega(y_1) = \omega(x_1)$ it follows that it must consist only of 
generators with $\omega(x) = 1$, i.e. of elements of $T$.

Therefore, $x_1$ and $y_1$ must differ by an integer linear combination of the elements of $T$, which is an element of $H_F$ as desired. 
\end{proof}

We now have a natural map from the geodesics with direction $F$ to the elements of the group $\ZZ^n/H_F$, given by taking a
geodesic to the coset of $H_F$ in $\ZZ^n$ containing all of its points. In fact, calling $\ZZ^n/H_F$ a group is somewhat
misleading; there is no natural group structure, since among other things there is no distinguished 
identity element. More
appropriate is to say that it is a set with a $\ZZ^n$-action on it.

Since for $G\subset F$, we have $H_G\subset H_F$, we might ask if there is a natural map from $\ZZ^n/H_G$ to $\ZZ^n/H_F$. Indeed
there is, and it is the one you would expect, the induced map between the sets as quotients of $\ZZ^n$; this map, on the
geodesics, is given by preserving the base and changing the direction. This clearly preserves the $\ZZ^n$-action; the preimage of
a point in $\ZZ^n/H_F$ is therefore just an $H_F$-orbit in $\ZZ^n/H_G$, whose compactification in $C(\ZZ^n, S)$ contains the
image point in question.

Furthermore, we can use this classification of points in $C(\ZZ^n, S)$ to prove that Cayley compactifications are actually compact.

\begin{thm}
The Cayley compactification $C(\ZZ^n, S)$ is compact for any choice of generating set $S$.
\end{thm}

\begin{proof}
We will show that any infinite sequence has a convergent subsequence.
Suppose we have an infinite sequence; because there are finitely many faces $F$, some infinite 
subsequence must be contained in an orbit $O_F$. So it suffices to show that any infinite sequence of geodesics with the same direction has 
a convergent subsequence.

For any such geodesic $\gamma$ with direction $F$, for sufficiently large $i$, a shortest path from 0 to $\gamma_i$ will consist of a 
shortest path from 0 to $\gamma_{i-1}$, followed by a path from $\gamma_{i-1}$ to $\gamma_i$. The second part of this consists only of 
elements of $F$, so to each $\gamma$ we can associate (not necessarily uniquely) a sum of generators not in $F$, which cannot be expressed 
as the sum of fewer generators. Suppose we have an infinite sequence of geodesics with direction $F$. Then for each generator $x\in S 
\setminus F$, either the coefficient of $x$ in the associated sum grows infinitely large, or there is some finite value of this coefficient 
which occurs infinitely often.

Thus, we can find a subsequence for which for each $x\in S$, either the coefficient of $x$ in the associated sums for these geodesics grows 
monotonically, or the coefficient of $x$ in these sums is constant. If $T$ is the set of generators for which the coefficient grows 
monotonically, then by a previous argument these geodesics converge to one with direction $F\cup T$, 
based at any point on any of these 
geodesics. Consequently, we have found a convergent subsequence, proving compactness.
\end{proof}

We have given a relatively complete description of the Cayley compactification using geometric arguments in the case of $\ZZ^n$. In the 
next section, we present a different method which will work in a more general case, although it does not give as complete results in the 
case of $\ZZ^n$. 

\section{Commutative algebra methods}\label{grobner}

In this section, we introduce commutative algebra methods, specifically
standard pair decompositions, that allow us to attack the question of Cayley compactifications in a more
general setting: in particular, when the group in question is not $\ZZ^n$, and when the generators have
different weights (which corresponds to the edges of the Cayley graph having different lengths depending
on which generator they arise from.) The language here also makes some of the proofs above simpler; we
present both versions of those proofs so as to betterdemonstrate the mechanics of both approaches. 

Given a Cayley graph $(G,S)$, we can form its Cayley ideal as follows.

\begin{defn}
Let $(G,S)$ be a Cayley graph. The Cayley ideal $I(G,S)$ is defined to be the kernel of the algebra homomorphism 
$k[S]\rightarrow k[G]$ sending each generator $s_i$ to the monomial $x^{s_i}$. Here, $k[S]$ is the free algebra with 
one generator for every element of $S$, while $k[G]$ is the group algebra with one generator $e_i$ for each element $g_i\in G$, and 
relations $e_ie_j = e_k$ whenever $g_i g_j = g_k$.
\end{defn}

The Cayley ideal is a binomial ideal, where each binomial represents two different paths from 0 to a
given point. We can therefore obtain a copy of the set of points by considering a basis of monomial
generators (monomials represent a product of generators, which corresponds to a point in $G$) for
$k[S]/I(S)$ as a $k$-vector space. A canonical way to do this is as follows: compute a Gr\"{o}bner basis for the ideal
$I(S)$, with respect to some cost vector $c$, where $c(s_i)$ represents the weight (i.e. length in the Cayley graph) of the generator.
Assuming $c$ is generic, this yields the initial monomial ideal $\text{in}(I(S))$. We recall the definition of a standard pair 
decomposition~\cite{stv}.

\begin{defn}
Let $J\subset k[t_1,\ldots,t_n]$ be a monomial ideal. A pair $(M, T)$, where $M$ is a monomial and $T$ is a subset of the 
variables, is defined to be a standard pair of $J$ if the following conditions hold.

(a) All monomials of the form $Mt_{i_1}^{e_1}t_{i_2}^{e_2}\cdots t_{i_k}^{e_k}$ are standard monomials (monomials not
contained in $J$), where each $t_{i_l}$ is in $T$,

(b) no variable appearing in the set $T$ appears in the monomial $M$, and

(c) $(M,T)$ is not properly contained in any other set with properties (a) and (b).
\end{defn}

It is clear from the definition that the union of the standard pairs of $J$ consists precisely of all the standard monomials of 
$J$.  We now claim that this standard pair decomposition of $\text{in}(I(S))$ gives us the elements of the
Cayley compactification $C(G,S)$ in a natural way.

\begin{lem}\label{stdmon}
The points of $G$ are in natural one-to-one correspondence with the standard monomials of $k[S]/I(S)$.
\end{lem}

\begin{proof}
A monomial in $k[S]$ represents a linear combination of generators; two linear combinations of generators represent the same point 
in $G$ if and only if the difference of the corresponding monomials is in $I(S)$. Because $S$ generates $G$,
each point can be 
represented as a linear combination of generators. But the set of standard monomials includes exactly one representative from each 
equivalence class of monomials in $k[S]/I(S)$, and hence exactly one representative for each point in $G$. Furthermore, the 
correspondence is natural, giving a shortest path from $0$ to that point in the weighted Cayley graph $C(G,S)$.
\end{proof}

Because of Lemma~\ref{stdmon}, by abuse of notation we will also consider the standard pairs as sets of points in $G$.

Now, suppose we have an infinite geodesic. There are only a finite number of standard pairs, so there must be some standard pair 
containing infinitely many points from that geodesic; the subset of the geodesic contained in this standard pair is a geodesic which 
is equivalent to the original one. Thus, when investigating the rest of the points in the Cayley compactification, we need only 
consider the geodesics whose points are entirely contained in a standard pair.

\begin{lem}
Suppose $\gamma = \{\gamma_0,\gamma_1,\ldots\}$ is a geodesic contained 
entirely in a standard pair $(M, \{s_1,\ldots,s_k\})$. Then 
$\gamma$ is equivalent to some geodesic $\eta$ with the property that for each $i$, $1\le i\le k$, one of the following is be the 
case:

(a) the exponent of $s_i$ in $\eta_j$ is constant, or

(b) the exponent of $s_i$ in $\eta_j$ grows monotonously without bound as $j$ gets large.
\end{lem}

\begin{proof}
Suppose $\gamma$ does not satisfy this property; suppose it violates it for $r$ different values of
$i$. We will induct on $r$. Take any $i$ for which the exponent of $s_i$ in $\gamma$ fails to obey
either condition. If the exponent of $s_i$ in $\gamma$ gets arbitrarily large, then we can find a
subsequence of $\gamma$ where this exponent is strictly increasing; this subsequence will then obey the
property for the same values $\gamma$ did, as well as for $i$, so it will violate the property for at
most $r-1$ different values of $i$. By induction, we are now done.

If not, then the exponent of $s_i$ in the infinite sequence $\gamma$ is bounded, so there exists some
value it takes infinitely often. Take the subsequence of $\gamma$ where it takes this value; this will
again obey the property for the same values $\gamma$ did, as well as for $i$, and we are again done by
induction.
\end{proof}

This further restricts the realm of geodesics. All infinite geodesics are equivalent to one contained in a standard pair
$(M,\{s_1,\ldots,s_k\})$, and furthermore all are equivalent to a geodesic $\gamma = \{\gamma_0,\gamma_1,\ldots\}$ where we have $\gamma_i
= N (t_1^{a_{1,i}} \cdots t_n^{a_{n,i}})$ where $N$ is some monomial in the standard pair not containing any of the $t_l$, the $t_l$ are a
subset of the $s_l$, and $a_{k,i} < a_{k,j}$ for $i<j$. Because of Theorem~\ref{basedir} from Section~\ref{boundary}, all such geodesics
with a point in common are equivalent. In fact, since we can prepend $N$ to any geodesic, all such geodesics with identical $N$'s are
equivalent; we call $N$ the \textbf{base} of the geodesic $\gamma$.

Therefore, using this standard pair decomposition, we can enumerate all of the points in the Cayley compactification, simply by considering
all possibilities for base and direction inside each standard pair. These possibilities are easy to enumerate; if $(M, \{t_1,\ldots,t_n\})$
is a standard pair, then the base-direction pairs are in one-to-one correspondence with expressions $Mt_1^{e_1}t_2^{e_2}\cdots t_n^{e_n}$,
where the $e_i$ range over $\NN\cup \{\infty\}$. 

Unfortunately, if the cost vector is not generic, it is not necessarily true 
that all of these geodesics are distinct. Indeed, we have the following example.

\textit{Example 4.} Consider the Cayley graph $(\ZZ^2, (\pm e_1, \pm e_2, \pm e_3, \pm e_4))$, where we
have $e_1 = (1,0)$, $e_2=(0,1)$, $e_3 = (1,1)$, and $e_4 = (1,-1)$. The Cayley compactification of this
is essentially a 4-gon, via the methods of Section~\ref{boundary}. However, the algebraic
approach yields a different answer. The standard pairs of the ideal $I(S)$ include $(1, \{e_3,e_4\})$ and
$(e_1, \{e_3,e_4\})$; we claim that the geodesics these code for are equivalent. Indeed, it is relatively
easy to check by hand that for a given $\phi_{x,y}$, along each of these geodesics its limit is equal to
the horizontal distance between $x$ and $y$.

While this method does not always work, it does generically give a complete answer.

\begin{thm}
Suppose that the cost vector is generic in the sense that $\text{in}(I(S))$ is a monomial ideal. Then no 
geodesics with different base or direction are equivalent.
\end{thm}

\begin{proof}
Suppose first that our two geodesics $\gamma_1 = (M,T)$ and $\gamma_2 = (N, U)$ have different 
directions; suppose without loss of generality that $t_1$ appears in $T$ but not $U$. 
Pick a point $x$ from $\gamma_1$ sufficiently far out that the exponent of $t_1$ in $x$ is higher than 
the exponent of $t_1$ in $N$; pick any point $y$ on $\gamma_2$. 

Now, we will consider the limit of the function $\phi_{x,y}$ along $\gamma_1$ and $\gamma_2$. Along 
$\gamma_1$ past $x$, the distance from a point $x_i$ to $x$ is just $\text{deg}(x_i)-\text{deg}(x)$, 
since both are standard monomials and $x_i$ is a multiple of $x$. 

Because $x_i$ is a standard monomial, the distance from $x_i$ to $y$ is at least 
$\text{deg}(x_i)-\text{deg}(y)$. Similarly, for $y_j$ sufficiently far out along $\gamma_2$, the 
distance from $y_j$ to $y$ is at $\text{deg}(y_j)-\text{deg}(y)$, while $d(y_j,x)$ 
is at least $\text{deg}(y_j)-\text{deg}(x)$. If $\gamma_1$ and $\gamma_2$ are equivalent, then we must 
for sufficiently large $i$ and $j$ have $d(x_i,x)-d(x_i,y)=d(y_i,x)-d(y_i,y)$; plugging in the above 
information, the only way this can happen is if we have $d(y_j,x)=\text{deg}(y_j)-\text{deg}(x)$ and 
$d(x_i,y)=\text{deg}(x_i)-\text{deg}(y)$.

This means that in the Cayley graph, there exists a path of length $\text{deg}(y_j)-\text{deg}(x)$
between $x$ and $y_j$. This path corresponds to a monomial $M_j$ of this degree. When we multiply this
monomial by $x$, we get another representation of degree $\text{deg}(y_j)$ for $y_j$. It cannot be the
same representation, since we chose $x$ such that the exponent of $t_1$ was higher in $x$ (and thus in
$xM_j$) than in $y_j$. So we have two distinct monomials of the same degree whose difference is in
$I(S)$, one of which is a standard monomial; this is impossible, since the cost vector is generic.

The argument in the case where the directions are the same and the bases are different is identical; as 
our point $x$, we simply pick whichever of $M$ and $N$ has a higher power of some element $s\in S$. 
This element will then have a higher power of $s$ than any element in the other geodesic, since by 
construction the variables appearing in the base and the direction to not overlap. The 
result then follows as in the first case.
\end{proof}

This allows us to effectively compute the Cayley compactification in the general case of an abelian group $G$ and a generic cost vector. 
The proof that the Cayley compactification is actually compact is easy with this approach; it is parallel to the proof we presented in 
Section~\ref{boundary}. 

\begin{thm}
For an abelian group $G$, and any set of weighted generators $S$, the Cayley compactification $C(G,S)$ is 
compact.
\end{thm}
\begin{proof}
Take any sequence of points in $C(G,S)$; we must show it has a convergent subsequence. Each point in 
the Cayley compactification belongs to a standard pair; there are only finitely many standard pairs, so 
some pair $(M, \{t_1,\ldots,t_k\})$ must contain an infinite subsequence. Now, recall that points of 
$C(G,S)$ correspond to expressions $Mt_1^{e_1}\cdots t_k^{e_k}$, where $e_i\in \NN\cup\{\infty\}$ for 
all $i$. All that remains is an elementary argument: either we can find a subsequence 
$\gamma_1,\gamma_2,\ldots$ where the value of $e_1$ increases monotonously without bound, or we can 
find one where the value of $e_1$ is a constant. We can now do this for each variable $e_i$ in turn to 
obtain a convergent subsequence (if some $e_i$ increases without bound then the limit has $e_i=\infty$; 
otherwise the limit has whatever value of $e_i$ all the elements do.) 
\end{proof}

\section{Conclusion} 

While we have used only undirected Cayley graphs for ease of notation, essentially
everything we have said in this paper holds for Cayley digraphs, where the generators yield directed
edges rather than edges, and we do not have $S=S^{-1}$. However, the case of nonabelian groups is
somewhat different; what can be said in this case? 

If we have a homomorphism between a group $G$ with generating set $S$ and a group $H$ with generating 
set $T$, and $S$ maps to $T$, there ought to be a relation of some kind between the Cayley 
compactifications. There are obstacles to this, however; geodesics in $G$ may not map to geodesics in 
$H$. For instance, if $(G,S)=(\ZZ, \{\pm 1\})$, and $(H,T) = (\ZZ, \{\pm 1, \pm 2\})$, geodesics in $G$ 
do 
not map to geodesics in $H$. If the map from $S$ to $T$ is surjective, there is some hope.

\section*{Acknowledgements}

The author was supported by an NSF graduate fellowship. I would also like to thank Marc Rieffel for bringing this problem to my attention,
and Bernd Sturmfels and David Eisenbud for suggesting fruitful directions of research.

\end{document}